\documentclass[12pt,english]{article}

\usepackage{graphicx,psfrag,epsfig,amsfonts,amssymb,amsmath,babel}

\newtheorem{Theorem}{\bf    Theorem}[section]
\newtheorem{Lemma}[Theorem]{\bf   Lemma}

\newtheorem{maintheorem}[Theorem]{\bf   Theorem}

\newtheorem{Definition}[Theorem]{\bf   Definition}
\newtheorem{maindefinition}[Theorem]{\bf  Definition}
\newtheorem{Example}[Theorem]{\bf   Example}

\newtheorem{Remark}[Theorem]{\bf   Remark}

\begin{document}

%


\title{{ On  Ilyashenko's Statistical Attractors}
\author{\rm Eleonora Catsigeras$^{\rm a}$ \thanks{   Instituto de Matem\'{a}tica y Estad\'{\i}stica \lq\lq Prof. Rafael Laguardia\rq\rq (IMERL). Facultad de Ingenier\'{\i}a -  Universidad de la Rep\'{u}blica. Address: Av. Herrera y Reissig 565, C.P. 11300,  Montevideo, Uruguay. Email: eleonora@fing.edu.uy }}}

\date{}


\maketitle

\begin{abstract}
We    define a minimal  $\alpha$-observability of Ilyashenko's statistical attractors. We prove   that   the
space is always full Lebesgue   decomposable into    pairwise
disjoint sets that are Lebesgue-bounded away from zero and included
in the basins of a finite family of  minimal $\alpha$-observable statistical attractors.
Among other examples, we analyze
   the   Bowen homeomorphisms
with   non robust topological heteroclinic cycles. We  prove the existence of
three types of statistical behaviours for   these examples.
\end{abstract}

\vspace{.5cm}



\section{Introduction}
The theory of the Ilyashenko's attractors was  principally developed  in \cite{3}, \cite{6},   \cite{4}, \cite{2}, \cite{1},
  \cite{5}, \cite{7}, \cite{8}. These attractors share the advantages of both Milnor's attractors \cite{Milnor}, \cite{Milnor2} and
  Pugh-Shub's ergodic attractors \cite{pughshub}. In fact, on the one hand, the sure existence of Milnor's attractors is inherited by Ilyashenko's attractors, since these latter exist for any continuous dynamical system on a compact Riemannian manifold (Theorem \ref{theoremExistenceOfAttractors}). On the other hand, the fine statistical description of Pugh-Shub's ergodic attractors is also inherited by any (minimal) Ilyashenko's attractor  $K$,  since any small neighbourhood of any point of $K$ must be visited in the future with a positive frequency  when time goes to infinite. As a counterpart, we notice that  in general  Pugh-Shub's ergodic attractors do not necessarily exist (Example \ref{exampleBowen}, Case C), while most points of Milnor's attractors may  be statistically irrelevant; namely,   small neighbourhoods  of them may receive     asymptotically zero-frequent visits in the future (Example \ref{examplehuyoung}).

   \begin{maindefinition}   \label{definitionMilnorLike} \em  \textsc{(Ilyas\-henko's Statistical Attractor)}

 Let $M$ be a compact Riemannian manifold. Let $f\colon M \to M$ be a \em continuous \em map. Let $K \subset M$ be a nonempty  compact   set. The set $K$ is a \em Ilyashenko's   statistical  attractor \em if:

  \noindent{\bf a) } the following set $B(K)$, which  is called  \em basin of statistical attraction \em of $K$  or in brief  \em basin \em of $K$, has positive Lebesgue measure:


 \begin{equation} \label{eqnB(K)} B(K) := \{x \in M: \ \ \lim_{n \rightarrow + \infty} \frac{1}{n} \, \#\{0 \leq j \leq n-1: \ \ \mbox{dist}(f^j(x), K) < \epsilon\} = 1 \ \ \forall  \epsilon >0\},  \end{equation}
 where $\# A $ denotes the cardinality of the finite set $A$;

\noindent{\bf b) } $K$ is minimal with respect to the basin $B(K)$, i.e.
$$K' \neq \emptyset \mbox{ is compact }, \ K' \subset K , \ B(K') = B(K) \ \  m-\mbox{a.e.} \ \ \Rightarrow \ \ K'= K, $$
where $m$ denotes the Lebesgue measure on $M$.
\end{maindefinition}

From Definition \ref{definitionMilnorLike},  the orbit of a point $x \in M$   belongs to the basin of statistical attraction
  $B(K)$ of $K$  if and only if  the asymptotic frequency
  of its visits in the future to any arbitrarily small neighbourhood of $K$ is 1.
  The future orbit does not need to remain near $K$ in all the
  instants of the large future.
  It is just required  that the frequency
  according to which one may find a  future iterate of $x$ far from  $K$,
  is asymptotically null. So, the attraction to $K$ is not necessarily topological, but statistical.

  Also notice that $K$ is a statistical attractor only if its basin $B(K)$ of statistical
  attraction has positive Lebesgue measure.
  Thus, a nonempty   basin $B(K)$ (as the compact support $K$ of any ergodic measure has) is not enough. It is immediate to check that if there exists an ergodic invariant measure $\mu$ that besides is physical, then the compact support $K$  of $\mu$ is an statistical attractor. Nevertheless, the converse is false, as Bowen example shows (cf. Example \ref{exampleBowen}).

\begin{Remark}
 \label{remarkf-invariance} \em
 Due to the continuity of $f$ and   to  Condition b) of minimality, it is deduced that any Ilyashenko's statistical attractor $K$ is $f$-invariant, i.e. $f^{-1}(K) = K$. In fact, take $y \in K$. After the minimality of $K$ with respect to its basin, for any open neighborhood $V \ni y$ there exists a positive-Lebesgue subset $A(y) \subset B(K)$ such that the frequency of visits to $V$ of the orbit of $x$ converges to 1, for all $x \in A(y)$. Conversely, if $y \in M$ is any point of the manifold such that the frequency of visits to $V$ of  the orbit  of $x$ converges to 1  for all $x$ in certain Lebesgue positive set $A(y) \subset B(K)$, then $y \in K$ (because  Equality (\ref{eqnB(K)}) and because the frequencies of visits to the neighborhoods of two disjoint  compact sets, cannot both converge to 1). Since $f$ is continuous, we deduce that
the frequency of visits of any orbit to a neighborhood $U$ of a  point $z = f(y)$ in the manifold, coincides with the frequency of visits of the same orbit to the neighborhood $V:= f^{-1}(U) \ni y$. Thus, we conclude that $y \in K$ if and only if $f(y) \in K$, proving that any atractor $K$ according to Definition \ref{definitionMilnorLike}, must be $f$-invariant. Nevertheless, we are interested to generalize the above definition \em for any (non necessarily continuous) Borel measurable map $f: M \mapsto M$. \em But  the argument  above does not work for non continuous $f$. On the one hand, if we imposed the $f$-invariance of $K$ in the definition, statistical attractors may not exist for non continuous $f$ (recall for instance the trivial example $f:[0,1] \mapsto [0,1]$ defined by $f(x) = x/2 $ if $x >0$ and $f(0)= 1$). On the other hand, if we did not impose the $f$-invariance to $K$, then statistical attractors would still   be characterized by means of   relevant probability measures  that generalize the concept of physical measures (Theorem \ref{theoremC}), but these measures would not necessarily be $f$-invariant.
\end{Remark}

\begin{maindefinition}   \label{definitionAlphaObserv}  \em \textsc{($\alpha$-Observability and $\alpha$-obs. Minimality)}

  For any $0 < \alpha \leq 1$, we say that a nonempty compact set $K$ satisfying condition a) of Definition \ref{definitionMilnorLike} is
  \em $\alpha$-observable \em  if $m(B(K)) \geq \alpha$, where $m$
  denotes the Lebesgue measure. We abbreviate this property by \em $\alpha$-obs.  \em
We say that   $K$ is \em minimal   $\alpha$-obs, \em if it
  is $\alpha$-obs. and   no proper compact subset of $K$  is also   $\alpha$-obs
  for the same value of $\alpha$.
 \end{maindefinition}

 In Remarks \ref{remarkMinimalidadComparacion0} and \ref{remarkMinimalidadComparacion} we discuss the relation between the $\alpha$-obs. minimality and the minimality condition (b) of Definition \ref{definitionMilnorLike}.

\begin{maindefinition}   \label{definitionMilnorLike2} \em \textsc{(Statistical attractor)}

 Let $M$ be a compact Riemannian manifold. Let $f\colon M \to M$ be a Borel measurable map. Let $K \subset M$ be a nonempty and compact set. We say that $K$ is a  \em   statistical attractor   \em  if it satisfies condition a) of Definition \ref{definitionMilnorLike} and besides:

 \noindent {\bf b') }  there exists $\alpha >0$ such that $K$ is minimal  $\alpha$-obs.

\end{maindefinition}

 \begin{Remark}
 \label{remarkMinimalidadComparacion0} \em
 Let us prove that Definitions \ref{definitionMilnorLike} and \ref{definitionMilnorLike2} are equivalent  if $f$ is continuous.

 First, any statistical attractor according to Definition \ref{definitionMilnorLike2}  is a Ilyashenko's statistical attractor according to Definition \ref{definitionMilnorLike}. In fact, let us see that Condition b) is   satisfied.  Take $K' \subset K$ nonempty and compact such that $m(B(K')) = m(B(K))  $. The condition $K' \subset K$ immediately implies $B(K') \subset B(K)$. Since $K$ is $\alpha$-obs. minimal, for some $\alpha >0$, then $m(B(K)) \geq \alpha$ and $m(B(K'')) < \alpha$ for any   compact nonempty  set $K''$ properly contained in   $K$. Since $m(B(K')) = m(B(K)) \geq \alpha$ we conclude that $K' = K$.

 Conversely, let us see that any statistical attractor $K$ satisfying Definition \ref{definitionMilnorLike} also satisfies Definition \ref{definitionMilnorLike2}.   In fact,    choose   $\alpha = m (B(K))$.
  Assume that $K$ satisfies Definition \ref{definitionMilnorLike} and that $K' \subset K$  is nonempty and compact such that $m(B(K')) \geq \alpha$. Since $K' \subset K$  then $B(K') \subset B(K)$. From $m(B(K')) \geq \alpha = m (B(K))$ we deduce $B(K')= B(K)$ Lebesgue a.e.  From Condition b) of Definition \ref{definitionMilnorLike} we deduce that $K'= K$. So, we have shown that $K$ is $\alpha$-obs. minimal.

 \end{Remark}

  \begin{Remark} \label{remarkMinimalidadComparacion} \em

We notice that
we are using the adjective \lq\lq $\alpha$-obs. minimal\rq\rq \ in the sense of a least set in the chain of inclusion of
  $\alpha$-observable nonempty compact sets, for a fixed value of $\alpha$ bounded away from zero.

  Let us see that the $\alpha$-obs. minimality for a \em previously specified value of $\alpha >0$ \em is indeed a restriction to the concept of statistical attractor. In other words, let us see that a statistical attractor  $K$ satisfying Definition \ref{definitionMilnorLike} (or equivalently Definition \ref{definitionMilnorLike2} for some $\alpha' >0$), is not necessarily $\alpha-$obs. minimal,  if one has a previously specified value $0 <\alpha < m(B(K))$.
  In fact, consider the Bowen example   (cf. Example \ref{exampleBowen}). It is  formed by an eye (homeomorhic to a two-dimensional compact disk) with two saddles $p$ and $q$ in its boundary and such that the interior of the eye is the basin $B_1$ of statistical attraction of the compact set $K_1 = \{p,q\}$. Now, add other Bowen example, i.e. other eye, whose two saddles are $q, r$, its interior $B_2$ is the basin of statistical attraction of $K_2 = \{q,r\}$, and such that the intersection of both eyes is only the saddle $q$. Assume for instance that $m(B_1) = m(B_2) = 1/2$. Then $K = \{p, q, r\}$ is a Ilyashenko's statistical attractor (satisfying Definition \ref{definitionMilnorLike} and equivalently Definition \ref{definitionMilnorLike2} with $\alpha= 1$), whose basin is $B= B_1 \cup B_2$. The statistical attractor $K$ is minimal \em with respect to its basin $B$  \em  and is also $1$-obs. minimal. But $K$ is not minimal in an absolute sense. In fact,  $K_1$ and  $K_2$ are also statistical attractors, whose basins are $B_1$ and $B_2$ respectively, and do also satisfy Definitions \ref{definitionMilnorLike} and \ref{definitionMilnorLike2}.   They are $1/2$-obs. minimal. In this example  $K_1$ and $K_2$ are    the unique $\alpha-$obs. minimal statistical attractors for any $0 < \alpha \leq 1/2$. They are proper compact subsets of $K$. So, the statistical attractor $K$ is not minimal (in the absolute sense)  among all the existing Ilyashenko's  attractors.

  The latter example is too simple, because  $K_1$ and $K_2$ are minimal in the absolute sense, among all the existing Ilyashenko's attractors. But in a general context, there may not exist   Ilyshenko's statistical attractors that are minimal among all the existing ones, unless a positive minimum $\alpha $ is previously specified for the Lebesgue measure of the basins of statistical attraction.
  In fact, consider in a three-dimensional setting the following example. According to a real parameter $\theta \in [0, \pi/2)$ immerse a two-dimensional Bowen example's eye with statistical attractor $K_{\theta} = \{0, p_{\theta}\}$,   with diameter going to zero when $\theta \rightarrow \pi/2$, and such that  the two-dimensional basin  $B(K_{\theta})$ of statistical attraction  (the interior  of each eye) is contained in a plane  forming angle $\theta$ with the horizontal plane. In this example, each Ilyashenko's statistical attractor, properly contains infinitely many  other attractors, even if each of them is minimal with respect to its basin. In spite of that, for each previously specified value $\alpha >0$, the space is still Lebesgue-a.e. decomposable into the basins of a \em finite number \em $N_{\alpha}$ of $\alpha$-obs. minimal statistical attractors (in this  example, if
  $\alpha$ goes to zero, then $N_{\alpha}$ goes to infinite).

  \end{Remark}

\newpage

  For the sake of completeness we include  the following definitions:

  \begin{Definition} \em
\label{definitionMilnorAttractor} \textsc{(Milnor's attractor)}

Let $f\colon M \to M$ be a Borel measurable map.
A compact set $K \subset M$ is a \em Milnor's attractor \em if the set $A(K) \subset M$ of all the initial states $x \in M$ such that the omega-limit set $\omega(x)$   is contained in $K$, has positive Lebesgue measure, and if $K$ is the minimal compact set that contains $\omega(x)$ for Lebesgue all the points   $x \in A(K)$.

We recall that $\omega(x)   $ is the compact nonempty set in $M$ composed by the limits of all the convergent subsequences of the orbit $\{f^n(x)\}_{n \in \mathbf{N}}$. We   call $A(K)$  \em the basin of topological attraction \em of $K$. We say that a Milnor's attractor is $\alpha$-observable if $m(A(K)) \geq \alpha$, where $m$ denotes the Lebesgue measure. We say that a Milnor's attractor $K$ is minimal $\alpha$-obs. if it is $\alpha$-obs. and no proper subset of $K$ is an $\alpha$-obs. Milnor's attractor.
 \end{Definition}

 \vspace{.2cm}

\noindent{\bf Notation: } Roughly abusing of the language, we will use the words \em statistical attractor \em   referring to   any nonempty compact set $K$ satisfying Condition a) of Definition \ref{definitionMilnorLike}, regardless whether   $K$ satisfies  Condition b) of minimality or not. If   besides $K$ is $\alpha$-obs. minimal for some $\alpha \in (0,1]$, then  $K$ also satisfies Condition b), and conversely (cf. Remark \ref{remarkMinimalidadComparacion0}). The rough use of the language will not produce a conflict with Definitions \ref{definitionMilnorLike} and \ref{definitionMilnorLike2} since we will always search for   such a compact set  $K$ that   besides is $\alpha$-obs. minimal for some $\alpha$.

  Since topological attraction implies statistical attraction, any $\alpha$-obs.  minimal Milnor's attractor     is an $\alpha$-observable (but maybe non minimal)  statistical attractor,
according to Definitions \ref{definitionAlphaObserv} and \ref{definitionMilnorLike2}.
But, as  we show    in Examples \ref{examplehuyoung} and
\ref{exampleBowen}, not all the minimal $\alpha$-obs. Milnor
attractors are   minimal $\alpha$-obs.  statistical attractors. Nevertheless,
as a corollary of Theorem \ref{TheoremD}   we prove the
following statement: The basin $A(K)$ of topological attraction of any $\alpha$-obs
Milnor's attractor $K$, is \em covered \em by the union   (up to a zero
Lebesgue measure set)   of the basins $B(K_i)$ of statistical attraction of a finite
family of minimal $\alpha_i$-obs  statistical attractors $K_i
\subset K$, for some adequate positive values of $\alpha_i$. In
the above result, the union of all the minimal statistical
attractors $K_i$ contained in $K$, is not necessarily   equal
to the Milnor's attractor $K$ (see Examples \ref{examplehuyoung}
and \ref{exampleBowen}). Therefore,   the statistical
attractors are thinner sets than the Milnor's attractors.

\vspace{.3cm}

 To state the following definition, we denote by ${\mathcal M}$ the space of all the Borel probability measures on $M$ endowed with the weak$^*$ topology. We denote by $\lim^*$   the limit in ${\mathcal M}$.
 
 \newpage

 \begin{Definition} \em
 \textsc{(empirical probabilities and the  limit set in  ${\mathcal M}$)} \label{definitionEmpirical}

 Let $x \in M$. The  \em sequence of empirical probabilities $\{\nu_n(x)\}_{n\geq 1}$ \em of $x  $ is defined by


\begin{equation}
\label{equation1}
   \nu_n(x) : = \frac{1}{n} \sum_{j= 0}^{n-1} \delta_{f^j(x)} \in {\mathcal M},
\end{equation}


\noindent being $\delta_y$ the Dirac delta supported on the point $y \in M$. In other words, $\nu_n(x)$ is the probability measure supported, and equally distributed, on the finite piece of the orbit of $x$ between the instant $0$ and the instant $n-1$.

The \em limit set in the space of probability measures \em ${\mathcal M}$ of the   orbit of $x \in M$ is:


 \begin{equation}{\mathcal L}^*(x) := \{ \mu \in {\mathcal M}: \ \exists \ n_i \rightarrow + \infty \mbox{ such that } {\lim_{i \rightarrow + \infty}}^{\!\!\!\!\!^*}\ {\nu_{n_i}(x)} = \mu\}   \label{equationL*x} . \end{equation}


Since ${\mathcal M}$ is compact and sequently compact, then   ${\mathcal L}^*(x)$ is
    nonempty  and compact  for all $x \in M$. We say that ${\mathcal L}^*(x)$
     describes \em the asymptotic statistics \em of the future orbit of $x$.

 \end{Definition}

\begin{Definition} \em
 \textsc{(srb measures and ergodic attractors)}
\label{definitionPhysical}

A probability Borel measure $\mu$ on $M$ is called \em SRB or physical \em if the set


$$B(\mu) := \{ x \in M: {\lim   _{n \rightarrow + \infty}}^{\!\!\!\!\! *} \ \nu_n(x) = \mu \} = \{x \in M: {\mathcal L}^*(x) = \{\mu\} \}$$


\noindent has positive Lebesgue measure.
We call $B(\mu)$ the \em basin of statistical attraction \em of $\mu$.

Abusing of the language (regardless whether   $\mu$ is   ergodic), we call \em ergodic attractor \em to
the compact support $K$ of an SRB measure $\mu$ (i.e. $K$ is
the minimal compact set such that $\mu(K) = 1$).

\end{Definition}

After Definition \ref{definitionPhysical}, any SRB  measure is $f$-invariant provided that $f$ is continuous. Nevertheless,  SRB  measures are not necessarily ergodic  (see for instance Bowen's Example \ref{exampleBowen}, Case (B) at the end of this paper).
We notice that the definition of SRB  measure also holds for any Borel-measurable map $f\colon M \to M$, but in this case $\mu$ is not necessarily $f$-invariant. For instance,   $f\colon [0,1] \to [0,1]$ defined by $f(0) = 1, \ f(x) = (1/2)x$ for all $x \neq 0$, has the SRB  measure $\delta_{0}$ (whose basin is  $ [0,1]$), which is not $f$-invariant.


 \subsection{Statement of the results} \label{sectionStatements}

  The contributions  of this paper to the theory of statistical attractors are:

  \noindent  {\bf a)} Definition \ref{definitionAlphaObserv} in which we introduce the concept of   minimal $\alpha$-observability of Ilyashenko's statistical attractors;

  \noindent  {\bf b)} Theorem  \ref{theoremExistenceOfAttractors}, which slightly strengthen the previously known result  of the existence of Ilyashenko's statistical attractors, by proving also their existence under the minimal $\alpha$-obs. condition for any previously specified value of $\alpha \in (0,1]$.

\noindent  {\bf c)} Theorem  \ref{theoremC}, which adds to  the previously known results derived from   the Krylov-\-Bogol\-yubov
proc\-edure, the relationship between SRB-like measures, defined in \cite{CE3},   and the statistical attractors of the system;

  \noindent  {\bf d)} Theorem \ref{TheoremD}, which constructs
  a natural decomposition of a Lebesgue-full subset of the  space
  into the  basins of a finite number of minimal observable statistical attractors;

  \noindent {\bf e)} the proof of the existence of three types of statistical behavior in a $C^0$-version of Bowen's diffeomorphisms (Example \ref{exampleBowen}).

\begin{maintheorem}
\label{theoremExistenceOfAttractors} \textsc{(Existence of $\alpha$-obs. statistical attractors)}  

Let $f\colon M \to M$ be   Borel-measurable. For all $0 < \alpha \leq 1$ there exist  minimal  $\alpha$-observable
statistical attractors.

   Moreover, if
$\alpha= 1$, then the minimal $\alpha$-obs.   statistical attractor is unique.
\end{maintheorem}

We prove   Theorem \ref{theoremExistenceOfAttractors} in Section \ref{SectionExistence}, where we also   prove Theorems
\ref{corollaryMinimalRestrictedM'},
\ref{corollaryMinimalAttractingB} and \ref{CorollaryH}, which are slight
generalizations of Theorem \ref{theoremExistenceOfAttractors} relative to   some  previously fixed invariant subsets of the    space.


 It is straightforward to check that the minimal compact support  of  an SRB measure, when
this measure exists, is an $\alpha$-obs. minimal statistical attractor
for some $\alpha>0$. The following Theorem
\ref{theoremC}  asserts a kind of converse statement:   any $\alpha$-obs minimal
statistical attractor is the minimal compact support of a set of SRB-like measures. The latter measures are obtained after applying the
Krylov-Bogolyubov procedure to the empirical probabilities
constructed in Equality (\ref{equation1}). The method
takes any weak$^*$ partial limit  of the time averages of  non necessarily
invariant  probabilities.

   \begin{maintheorem}
\label{theoremC} \textsc{(Characterization of $\alpha$-obs statistical attractors)}

If $K$ is a   statistical attractor  and  if $ B(K)  $  is its  basin, then there exists a unique non empty weak$^*$-compact set ${\mathcal O}_f(K)$ of  probability measures  \em (which we call SRB-like measures) \em  such that:


 \em (a) \em For Lebesgue almost all the initial states $x \in B(K)$, and for all the   convergent subsequences  of  the  empirical distributions $ \nu_n(x)  : = (1/n) \, \sum_{j= 0}^{n-1} \delta_{f^j(x)}$, their weak$^*$-limits are probability measures contained in $ {\mathcal O}_f(K)$.


 \em (b) \em ${\mathcal O}_f(K)$ is the minimum weak$^*$-compact set of probability measures satisfying \em (a).


   (c) \em  $\mu(K) = 1$ for all $\mu \in {\mathcal O}_{f}(K)$.

  \em  (d) \em  If besides $K$ is minimal $\alpha$-obs. for some $0 < \alpha \leq 1$, then $K$ is the  minimal compact set in $M$ such that $\mu(K) = 1 \ \forall \ \mu \in {\mathcal O}_f(K)$.

\end{maintheorem}

We   prove Theorem \ref{theoremC} in Section \ref{sectionprooftheoremC}. In the proof we use
  the   SRB-like
measures, defined in \cite{CE3},  which  we restate in Definition \ref{definitionobservable}.  In Section \ref{sectionprooftheoremC} we also prove the converse of part (d) of Theorem \ref{theoremC}. In fact, we show that the minimal compact support of all the measures $\mu$ in the weak$^*$-compact set of   SRB-like measures,  is a minimal $\alpha$-obs. statistical attractor for some $0 < \alpha \leq 1$ (see Theorem \ref{CorollaryH}).

Theorem \ref{theoremC} states that for all $0<\alpha\leq 1$,
any $\alpha$-observable statistical attractor $K$ is provided
with a minimal weak$^*$-compact subset ${\mathcal O}_f(K)$ of
probability measures  with two remarkable properties:

(1) It is a set of $f$-invariant measures which has,
with respect to the attractor $K$, a  \lq \lq physical\rq\rq \ role
as SRB measures have, when they exist,  with respect to the
  ergodic attractors. In fact, after the statement (a)
of Theorem \ref{theoremC}, the invariant measures in ${\mathcal
O}_f$   completely describe the asymptotic statistics of the time series for
Lebesgue-almost all the orbits attracted by $K$.

(2) It is the minimal compact set of probability measures that
completely describes the asymptotic statistics, as
part (b) of Theorem \ref{theoremC} states. Therefore, the statistical
attractors   are the optimal choice, among the compact invariant
sets in the ambient manifold $M$,  if one aims to describe
the Lebesgue-full  asymptotic statistics of the system.

\vspace{.2cm}

The following Theorem \ref{TheoremD}  states  the existence and finitude of a
  decomposition of the space, up to  a  zero-Lebesgue subset,
 into a  family of   sets, each one contained in the basin  of
attraction of a  statistical attractor satisfying a minimally observable condition.

\begin{maintheorem} \label{TheoremD} \textsc{(Finite decomposition into statistical attractors)}

Let $0 < \alpha \leq 1$ be fixed. Let $m$ denote the Lebesgue probability measure.

There exists a finite family $\{K_i\}_{1 \leq i \leq p}$ of  $\alpha_i $-obs statistical attractors $ K_i $  with basins $B(K_i)$ such that:

 \em (a) \em    $\bigcup_{i=1}^p B(K_i)$ covers $m-$almost all $M$.

  \em (b) \em  $\alpha_i  = \alpha$ for all the values of $i \in \{1, \ldots, p\}$
  except at most one, say $  i_0 $, for which  $0 <\alpha_{i_0} = m(B(K_{i_0})) < \alpha.$   \em
(Therefore
   $m(B(K_i)) \geq \alpha \ \ \forall \ i \in \{1, \ldots, p\}      \mbox{  such that }   i \neq i_0.$)

     (c) \em For all $1 \leq i \leq p$ the statistical attractor $K_i$ is $\alpha_i$-obs. minimal for $f$ restricted to $M \setminus \cup_{j= 1}^{i-1} B(K_j)$. \em  (We denote $\cup_{j= 1}^{0} \cdot = \emptyset$.)
\end{maintheorem}

We   prove Theorem \ref{TheoremD} in Section
\ref{proofTheoremD}.  The  proof is rather natural:   roughly speaking,  one can take away  minimal
 observable sets (together with what they attract), one by one.

We notice that the statistical attractors $K_i$ of the
decomposition in Theorem \ref{TheoremD}  are not necessarily pairwise disjoint.
If all the statistical attractors $K_i$ are   mutually disjoint,
then any pair of them   would be at positive distance
(since they are compact sets), and so their
basins would be also pairwise disjoint.
If the additional assumption of mutually disjointness of the attractors holds,
Theorem \ref{TheoremD} asserts that
the basins $B(K_i)$ of the finitely many statistical attractors
$K_i$ would form a partition of   Lebesgue-a.e. the  space.
Anyway,  if the disjointness condition does not hold,
the basins of attractions of the finite number of
statistical attractors cover Lebesgue-a.e. the space,
and are, one by one,  Lebesgue-bounded away from zero.

To end this section, we deduce an immediate corollary of Theorem
\ref{TheoremD}, which shows  that the statistical attractors
are thinner  than   Milnor's attractors: namely, each $\alpha$-obs. minimal Milnor's attractor contains the union of a finite number of   statistical attractors.

\vspace{.3cm}

\noindent{\bf Corollary of Theorem
\ref{TheoremD}}

\em Let $0 < \alpha \leq 1$, and let $K$ be an $\alpha$-obs. minimal Milnor's attractor with basin $A(K)$.
There exists a finite    number of statistical attractors $K_1, \ldots, K_p$ contained
in $K$  that satisfy the conditions \em (a), (b) \em and \em (c) \em of Theorem \em \ref{TheoremD} \em for the set $A(K)$ instead of $M$. \em


This corollary is   immediate after
   Theorem \ref{TheoremD}. In fact,    along the proof of Theorem \ref{TheoremD}
   one does not use  the manifold
   structure of the ambient space $M$  for any purpose except   to define its Lebesgue measure $m$.
  Therefore, to prove the corollary it is enough to put   $f|_{A(K)}\colon A(K)
\to A(K)$ in the role of $f\colon M \to M$ and $m|_{A(K)}$ in
the role of $m$, where $m|_{A(K)}:= m(B \cap A(K))$
for any Borel set $B \subset M$.


\section{Proof of the existence   of minimal $\alpha$-obs. statistical attractors} \label{SectionExistence}

In this section we prove Theorem
\ref{theoremExistenceOfAttractors}. We   also introduce  some definitions which  impose  additional
minimal conditions to the statistical attractors
(Definitions \ref{definitionMinimalRestrictedM'} and
\ref{definitionMinimalAttractingB}). At the end of this section
we strengthen Theorem \ref{theoremExistenceOfAttractors}, proving the existence of statistical attractors satisfying those additional
conditions (Theorems \ref{corollaryMinimalRestrictedM'} and
\ref{corollaryMinimalAttractingB}).

%

\vspace{.3cm}
\noindent  {\bf Proof  of Theorem
\ref{theoremExistenceOfAttractors}}

{\em Proof: }
Let us fix $0 <\alpha \leq 1$. Consider the family $\aleph_{\alpha}$ of all the $\alpha$-obs statistical attractors (non necessarily minimal). The family $\aleph_{\alpha}$ is nonempty  since it trivially contains the manifold $M$.

Define in $\aleph_{\alpha}$ the partial order $K_1 \leq K_2$ if $K_1 \subset K_2$. Since the
attractors are all non empty compact sets, any chain $\{K_a\}_{a \in A  } \subset \aleph_{\alpha}$ (i.e. any totally ordered subset of $\aleph_{\alpha}$) has a non empty compact intersection:
$K := \bigcap_{a \in A} K_a.$   Let us prove that $K \in \aleph_{\alpha}$.

We have
to prove that   $m(B(K))\geq \alpha$,
where $m$ is the Lebesgue measure and $B(K)$ is the basin of statistical attraction of $K$,
constructed  in Definition \ref{definitionMilnorLike}.

For any $\epsilon >0$  and for any nonempty compact set $H \subset M$, define


\begin{equation}
\label{equationBsubEpsilonDeH}
B_{\epsilon}(H) := \{x \in M: \lim_{n \rightarrow + \infty} \omega_{n,H, \epsilon}(x) = 1\}, \ \ \mbox{ where }\end{equation}


$$\omega_{n,H, \epsilon} (x) := \frac{1}{n} \,   \, \#\{0 \leq j \leq n-1: \ \mbox{dist}(f^j(x), H) < \epsilon\} \leq 1.$$


\noindent It is standard to check that $B_{\epsilon'} (K) \subset B_{\epsilon}(K)$ if $0 < \epsilon' < \epsilon$. Therefore,   $$B(K) = \bigcap _{\epsilon>0} B_{\epsilon}(K) = \bigcap_{n \geq 1} B_{1/n}(K),$$  and thus $$m (B(K)) = \lim_{n \rightarrow + \infty} m(B_{1/n}(K)) = \lim _{\epsilon \rightarrow 0^+} m(B_{\epsilon}(K)).$$ So, to deduce that $m(B(K)) \geq \alpha$ it is enough to prove that $m(B_{\epsilon}(K)) \geq \alpha$ for all $\epsilon >0$.

Fix $\epsilon >0$. We assert that there exists $a \in A$ such that $\mbox{dist}(y, K) < \epsilon$ for all $y \in K_a$. Arguing by contradiction, if the intersection $K_a \bigcap \{y \in M: \mbox{dist}(y, K) \geq \epsilon\} $ were nonempty for all $a \in A$, since $\{K_a\}_{a \in A}$ is totally ordered, the property of finite intersections of compact sets would imply that $ \bigcap _{a \in A} K_a \bigcap \{y \in M: \mbox{dist}(y, K) \geq \epsilon\} \neq \emptyset $, contradicting the construction of $K= \bigcap _{a \in A} K_a$ and proving the assertion.

Using the triangle property for the value of $a \in A$ satisfying the above assertion,   we deduce $\omega_{n, K_a, \epsilon} (x) \leq \omega_{n, K, 2 \epsilon} (x)$ for all $x \in M$ and for all $n \in \mathbb{N}$. Therefore $B_{\epsilon}(K_{a}) \subset B_{2 \epsilon}(K)$. Since $K_{a} \in \aleph_{\alpha}$, we obtain $\alpha \leq m(B(K_a)) \leq m(B_{\epsilon}(K_a)) \leq m (B_{2 \epsilon}(K))$, as wanted. We have proved that $K \in \aleph_{\alpha}$ and so, any chain in $\aleph_{\alpha}$ has a minimal element.  Applying Zorn Lemma we deduce that there exist minimal elements in $\aleph_{\alpha}$. This means that there exist
$\alpha$-obs statistical attractors $K \subset M$, that do not
contain  proper subsets that  are also $\alpha$-obs statistical
attractors. So, the existence of minimal $\alpha$-obs statistical attractors is proved for any
previously specified value of $\alpha \in (0, 1]$.

To end the proof of Theorem \ref{theoremExistenceOfAttractors} it is left to show that the
minimal $1$-obs. statistical attractor is unique. In fact, consider $K_1$ and $K_2$, minimal $1$-obs. statistical attractors, namely their basins $B(K_1)$ and $B(K_2)$ have full Lebesgue measures: $m (B(K_1)) = m (B(K_2)) = 1$. Therefore, $m(B(K_1) \cap B(K_2)) = 1$. Take $x \in B(K_1) \cap B(K_2)$. Thus $1= \lim_{n \rightarrow + \infty} \omega_{n, K_i, \epsilon}(x)  $  for all $\epsilon >0$, for $i= 1, 2$.
Denote


$$\omega_{n, K_1, K_2, \epsilon}(x) = \frac{1}{n} \#\{0 \leq j \leq n-1:  \mbox{dist}(f^j(x), K_i) < \epsilon \ \mbox{ for } i= 1, 2\}.$$


\noindent We have
$\omega_{n, K_1, \epsilon}(x) \leq \omega_{n, K_1, K_2, \epsilon}(x) + \big (1 - \omega_{n, K_2, \epsilon}(x) \big).$
Therefore, $\lim_{n \rightarrow + \infty} \omega_{n, K_1, K_2, \epsilon}(x) = 1  $ for Lebesgue almost all $x \in M$.  Besides, if $\mbox{dist}(y, K_i) < \epsilon$ for $i= 1, 2$, then
$\mbox{dist}\big (y, K_2 \cap\{z \in M: \mbox{dist}(K_1, z) \leq 2 \epsilon\} \big ) < \epsilon$. Thus,


$$\omega_{n, K_1, K_2, \epsilon}(x) \leq \omega_{n, K_2 \cap \{z \in M: \mbox{ dist} (z, K_1) \leq 2 \epsilon\}, \epsilon}(x) \leq \omega_{n, K_2 \cap \{z \in M: \mbox{ dist} (z, K_1) \leq 2 \epsilon_0\}, \epsilon}(x) $$


\noindent for all $ 0 < \epsilon \leq \epsilon_0.$ We deduce that, for each fixed $\epsilon_0 >0$, the compact set $K_2 \cap \{z \in M: \mbox{ dist} (z, K_1) \leq 2 \epsilon_0\} $ is an $1$-obs. statistical attractor. Since $K_2$ is minimal with such a property, we obtain
$K_2 \subset \{z \in M: \mbox{ dist} (z, K_1) \leq 2 \epsilon_0\}  $ for all $\epsilon_0 >0$. Therefore, $K_2 \subset K_1$. Arguing symmetrically, $K_1 \subset K_2$, and thus  $K_1 = K_2$ ending the proof.    \hfill $\Box$


\begin{Definition}   \label{definitionMinimalRestrictedM'} \em
Let $0 \leq \alpha \leq 1$ and let $M' \subset M$ be a Borel set such that $M' \subset f^{-1}(M')$ and $m(M') \geq \alpha$. We say that a nonempty, compact and $f$-invariant set $K \subset M$   is an \em $\alpha$-obs  statistical attractor restricted to $M'$, \em or \em for $f|_M$, \em if its basin $B(K) $, as defined in \ref{definitionMilnorLike}, satisfies:


\begin{equation}
\label{equationMinimalRestrictedM'}
m(B(K) \bigcap M') \geq \alpha.\end{equation}


We say that an $\alpha$-obs statistical attractor $K$ is \em minimal restricted to $M'$, \em or \em for $f|_M$,  \em
if it satisfies the inequality (\ref{equationMinimalRestrictedM'}) and
has not proper, nonempty and compact subsets that satisfy it.
\end{Definition}
\begin{Definition} \em
\label{definitionMinimalAttractingB}
Let $B \subset M$ be a Borel set such that $B \subset f^{-1}(B)$ and $m(B) \geq \alpha >0$. We say that a nonempty compact and $f$-invariant set $K \subset M$ is a  \em  statistical attractor attracting $B$ \em if its basin of attraction $B(K)$, as defined in \ref{definitionMilnorLike}, satisfies:


\begin{equation}
\label{equationMinimalAttractingB}
B(K) \supset B \ m-\mbox{a.e.}  \ \ \mbox{ In other words, } \ m(B \setminus B(K)) = 0.
\end{equation}


 Since $m(B) \geq \alpha $  any statistical attractor attracting $B$ is $\alpha-$obs.

 We say that a statistical attractor is \em minimal attracting $B$ \em if it satisfies the condition (\ref{equationMinimalAttractingB}) and has not proper, nonempty and compact subsets that satisfy it.

\end{Definition}

It is standard to check that an $\alpha$-obs minimal statistical attractor $K$ is also $\alpha$-obs minimal restricted to its basin, and  minimal attracting its basin.
\begin{maintheorem} \label{corollaryMinimalRestrictedM'}  
Let $M' \subset M$ be a Borel set such that $M' \subset f^{-1}(M')$ and $m(M') \geq \alpha >0$. Then, there exists an $\alpha$-obs statistical attractor that is minimal restricted to $M'$, according to  Definition \em \ref{definitionMinimalRestrictedM'}.
\end{maintheorem}

{\em Proof: }
Repeat the   proof of Theorem
\ref{theoremExistenceOfAttractors},   using $M'$ in the role of
$M$, $B(K) \cap M'$ in the role of $B(K)$ and $B_{\epsilon} (H)
\cap M'$ in the role of $B_{\epsilon}(H)$.
\hfill $\Box$

\begin{maintheorem}   \label{corollaryMinimalAttractingB}
Let $B  \subset M$ be a Borel set such that $B \subset f^{-1}(B)$ and $m(B)  >0$. Then, there exists a statistical attractor that is minimal attracting $B$, according to  Definition \em \ref{definitionMinimalAttractingB}.
\end{maintheorem}

{\em Proof: }
Repeat the proof of Theorem
\ref{theoremExistenceOfAttractors},   defining the family
$\aleph _{B}$  (instead of $\aleph_{\alpha}$)  of all the
statistical attractors $K \subset M$ such that $B(K) \supset B
$ \ $m-$a.e.   \hfill $\Box$


\section{Proof of the probabilistic characterization of statistical attractors} \label{SectionCharacterization} \label{sectionprooftheoremC}


In this section we prove Theorem \ref{theoremC}. To do that, we first
revisit the   definition of   SRB-like measures, taken from \cite{CE3}.
Let us fix a metric $\mbox{dist}^*$ inducing the weak$^*$ topology in the space ${\mathcal M}$ of all the Borel probability measures on $M$.

 \begin{Definition}
 \textsc{(SRB-like measures)} \em
 \label{definitionobservable}

 Let $B \subset M$ be a forward invariant set (i.e. $B \subset f^{-1}(B)$) that has positive Lebesgue measure. We say that a  probability measure $\mu$ is \em SRB-like or physical-like \em  for $f|_{B}$, if for all $\epsilon >0$ the following set $B_{\epsilon}(\mu) \subset B$ has positive Lebesgue measure:


 $$B_{\epsilon}(\mu) := \{x \in B: \ \mbox{dist}^*({\mathcal L}^*(x) , \mu) < \epsilon\},$$


 \noindent where ${\mathcal L}^*(x)$ is the nonempty weak$^*$-compact set defined in (\ref{equationL*x}).

 We call $B_{\epsilon}(\mu)$  \em the basin of $\epsilon$-weak statistical attraction \em of the probability $\mu$. We
  denote by ${\mathcal O}_{f| {B}}$    the set of all the SRB-like measures $\mu$ for $f|_{B}$.
To justify the name \lq\lq SRB-like measures\rq\rq,
compare Definition \ref{definitionobservable} with Definition \ref{definitionPhysical}.

If $f$ is continuous, then all the measures in  ${\mathcal O}_{f|B}$  are $f$-invariant. In other words ${\mathcal O}_{f|B}   \subset{\mathcal M}_f$. In fact, ${\mathcal L}^*(x) \subset {\mathcal M}_f$  for all $x \in M $, so any $\mu \in {\mathcal O}_{f|B}$ belongs to the weak$^*$-closure of ${\mathcal M}_f$. But if $f$ is continuous, then ${\mathcal M}_f$ is weak$^*$-closed. Thus, $\mu \in {\mathcal M}_f$  as wanted.
 \end{Definition}
 \noindent The lemmas \ref{lemma1} and \ref{lemma2} below, are reformulations of   results communicated
 in \cite{CE3}.
\begin{Lemma} \label{lemma1}
${\mathcal O}_{f|_B}$ is weak$^*$-compact and nonempty.
\end{Lemma}
{\em Proof: }
It is immediate that ${\mathcal O}_{f| B}$ is weak$^*$-compact, because it is closed in the   space ${\mathcal M}$, which is a compact metric space for any metric inducing its weak$^*$ topology.
Let us prove that it is nonempty. Assume by contradiction that no measure in $\mathcal M$
is SRB-like. Then for all $\mu \in \mathcal M$ there exists $\epsilon >0$ such
that $m(B_{\epsilon}(\mu)) = 0$, where $m$ denotes the Lebesgue   measure on $M$.
Since ${\mathcal M}$ is compact, there exists a finite covering of ${\mathcal M}$ with
balls $\{{\mathcal B}_i\}_{i= 1, \ldots, s}$ of radii $\epsilon_i$, centred at $\mu_i$ and such that $m(B_{\epsilon_i}(\mu_i))= 0$ for all $i = 1, \ldots, s$.   Since $m(\bigcup_{i=1}^s B_{\epsilon_i}(\mu_i) )= 0$ and  $ \bigcup_{i=1}^s B_{\epsilon_i}(\mu_i) \supset \{x \in B(K): \ \ {\mathcal L}^*(x)  \bigcap {\mathcal M} \neq \emptyset\}$,
we conclude that for Lebesgue almost all $x \in B(K)$ the limit set ${\mathcal L}^*(x)  $
(which by definition is always contained in the space ${\mathcal M}$), is empty.
This is a contradiction since the space ${\mathcal M}$ is sequentially
compact when endowed with the weak$^*$ topology, and thus,
${\mathcal L}^*(x) \neq \emptyset$ for all $x \in B(K)$.
  \hfill $\Box$

\begin{Lemma} \label{lemma2}
The set ${\mathcal O}_{f|B}$ is the minimum weak$^*$ compact set
in the space ${\mathcal M}$ of Borel probability measures such that ${\mathcal L}^*(x)
\subset {\mathcal O}_{f|B}$ for Lebesgue almost all $x \in B$.
\end{Lemma}
{\em Proof: }
Let us first prove that for $m$-almost all $x \in B$ the limit set ${\mathcal L}^*(x)$ is contained in ${\mathcal O}_{f|B}$. Assume by contradiction that the set of such points $x$ has $m$-measure smaller than   $  m(B)$. Then  $\lim _{\epsilon \rightarrow 0 } m(A_{\epsilon})  < m(B) ,$ where


$$A_{\epsilon}  := \{x \in B: \ \ \max  \{ \mbox{dist}^*(\nu, \mu):  \   \nu \in {\mathcal L}^*(x), \, \mu \in {\mathcal O}_{f|B} \} \ < \ \epsilon \}.$$


Then, for some $\epsilon_0 >0$ small enough
$m(B \setminus A_{\epsilon_0 })> 0  $.
In other words, for a Lebesgue positive set of points $x \in B$, the limit set ${\mathcal L}^*(x)$ intersects   the weak$^*$-compact set ${\mathcal K} := \{ \mu \in {\mathcal M} :  \ \mbox{dist}^*(\mu,   {\mathcal O}_{f|B}) \geq \epsilon_0\}$.  Therefore, at least one of the measures $\mu \in {\mathcal K}$ satisfies $m(B_{\epsilon}(\mu)) >0$ for all $0 < \epsilon \leq \epsilon_0$, where


 $$B_{\epsilon}(\mu):= \{x \in B: \ \ \mbox{dist}^*({\mathcal L}^*(x), \mu) < \epsilon \}.$$


 In fact, if the latter assertion were not true, we would cover ${\mathcal K}$ with a finite number of balls $\{{\mathcal B}_i\}_{i = 1, \ldots, s}$ such that for Lebesgue almost all point $x \in B$, $\ {\mathcal L}^* (x) \bigcap {\mathcal B}_i = \emptyset $ for all $i = 1, \ldots, s$. Thus ${\mathcal L}^* (x) \bigcap {\mathcal K} = \emptyset$ for Lebesgue almost all $x \in B$, contradicting the construction of the set ${\mathcal K}$.

Thus, there exists $\mu \in {\mathcal K} $ such that
$m(B_{\epsilon}(\mu)) >0$ for all $0 < \epsilon \leq \epsilon_0$
Then, after Definition \ref{definitionobservable} the probability measure $\mu$
is SRB-like for $f|B$. Therefore ${\mathcal K} \bigcap {\mathcal O}_{f|B} \neq \emptyset$,
contradicting the construction of the   compact set ${\mathcal K}$.
This ends the proof of the first assertion: \ for $m$-almost all $x \in B$, ${\mathcal L}^*(x) \subset {\mathcal O}_{f|B}.$

Second, let us prove that ${\mathcal O}_{f|B}$ is   minimal among all the weak$^*$ compact sets containing ${\mathcal L}^*(x)$ for Lebesgue almost all $x \in B$. In fact, if $\emptyset \neq {\mathcal K} \subset {\mathcal O}_{f|B} $, and $\mathcal K$ is compact, any measure $\mu \in {\mathcal O}_{f|B}\setminus {\mathcal K}$ is at positive distance, say $\epsilon >0$ (depending on $\mu$), from ${\mathcal K}$. After Definition \ref{definitionobservable} there exists a $m$-positive set of points $x \in B$ such that
$\mbox{dist}^*({\mathcal L}^*(x) , \mu) < \epsilon$. Therefore ${\mathcal L}^*(x) \not \subset {\mathcal K}$ for those points $x$. We conclude that ${\mathcal O}_{f|B}$ has not a proper, nonempty and compact subset  ${\mathcal K}$  containing ${\mathcal L}^*(x)$ for Lebesgue almost all $x \in B$. This ends the proof that ${\mathcal O} _{f|B}$ is minimal with such a property.
  \hfill $\Box$

\begin{Lemma}
\label{lemma3}  
If $K \subset M$ is a compact set such that $\mu(K)= 1$ for all $\mu \in \mathcal O_{f|B}$,   then $K $ is a statistical attractor whose basin $B(K)$ contains $B$ Lebesgue a.e.

\end{Lemma}

{\em Proof: }
Fix $\epsilon >0$ and choose any continuous function $\psi \in C^0(M, [0,1])$ such that $\psi|K = 1$ and $\psi(y) = 0 $ for all $y \in M$ such that $\mbox{dist}(y, K) \geq \epsilon$.
Choose and fix   $x \in B$, and a sequence of natural numbers $n_i \rightarrow + \infty$ such that
the following limits exist:


$$L = \lim _{i \rightarrow + \infty} \frac{1}{n_i }  \,  \# \{0 \leq j \leq n_i-1: \ \mbox{dist}(f^j(x), K) < \epsilon\}. $$


 $$\mu = {\lim   _{i \rightarrow + \infty}}^{\!\!\!\!*} (1/n_i) \sum_{j=0}^{n_i-1} \delta_{f^j(x)}  \in {\mathcal L}^*(x).$$


   On the one hand, ${\mathcal L}^*(x) \subset {\mathcal O}_{f|B}$ for $m$-a.e. $x \in B$, due to Lemma \ref{lemma2}. Besides, by hypothesis,    $\mu(K)=1$ for all   $\mu \in {\mathcal O}_{f|B}$. Therefore, the expected value of $\psi$ respect to the probability $\mu$ is equal to 1. In fact:
$1 \geq \int \psi \, d \mu \geq  \int_K \psi \, d \mu = \mu (K) = 1. $
On the other hand, since $\psi$ is continuous, the weak$^*$-limit in the space of probability measures implies:


$$1 = \int \psi \, d \mu = \lim _{i \rightarrow + \infty} \int \psi \, d \left ( \frac{1}{n_i} \sum _{j= 0}^{n_i-1} \delta_{f^j(x)}  \right )=   \lim_{i \rightarrow + \infty} \frac{1}{n_i} \sum_{j= 0}^{n_i-1} \psi (f^j(x)), $$


\noindent and, by construction of   $\psi$:


  $$\ \ \frac{1}{n } \, \sum_{j= 0}^{n -1} \psi (f^j(x)) \leq ({1}/{n}) \# \{0 \leq j \leq n-1: \ \mbox{dist}(f^j(x), K) < \epsilon\} \leq 1.$$

\noindent Then,
$1= \lim_{n \rightarrow + \infty} \frac{1}{n  }  \,  \# \{0 \leq j \leq n -1: \ \mbox{dist}(f^j(x), K) < \epsilon\} $ for $      m\mbox{-a.e.} \ x \in B.$
We deduce that $x \in B (K)$ for Lebesgue almost all  $x \in
B$, and so $K$ is a statistical attractor whose basin contains
Lebesgue a.e. $B$.    \hfill $\Box$

 \vspace{.4cm}

 \noindent {\bf End of the proof of Theorem   \ref{theoremC}.}


{\em Proof: }
    Consider the basin of attraction $B(K)$ of a given minimal $\alpha$-obs. statistical attractor $K$. By hypothesis $m(B(K)) \geq \alpha >0$. It is straightforward to check that if $x \in B(K)$ then $f(x) \in B(K)$ (even if $f$ is only a measurable map that is not continuous). Then, we can apply Definition \ref{definitionobservable}, and consider the set ${\mathcal O}_{f|B(K)}$ of all the SRB-like measures for $f|_{B(K)}$. After Lemmas \ref{lemma1} and \ref{lemma2} (denoting ${\mathcal O}_f(K)$ to  ${\mathcal O}_{f|B(K)}$), Assertions (a) and (b) of Theorem \ref{theoremC} are proved.

  Now, let us prove Assertion  (c).
   We shall prove that $\mu(K)= 1$ for all $\mu \in {\mathcal O}_f(K)$. Fix $\mu \in {\mathcal O}_f(K)$, choose an arbitrarily small $\epsilon >0$ and denote


  $$V_{\epsilon} = \{x \in M: \; \mbox{dist}(x, K) < \epsilon\}.$$


   Construct a  continuous real function $\psi  \in C^0(M, [0,1])$ such that $\psi|_K = 1$ and $\psi (x) = 0$ if $x \not \in V_{\epsilon}$. After the continuity of $\psi$ there exists $0 <\epsilon' < \epsilon$ such that $\psi(x) > 1 - \epsilon \ \ \forall \ x \in V_{\epsilon'}(K). $  Let us compute   the expected value of $\psi$ respect to the probability $\mu$:


 \begin{equation}
  \label{equation7}
  \int \psi \, d \mu  = \int_{V_{\epsilon}}  \psi \, d \mu \leq  \mu (V_{\epsilon}).\end{equation}


 Recall Equality (\ref{equationL*x}) which defines ${\mathcal L}^*(x)$ for all $x \in M$ and   Definition \ref{definitionMilnorLike} of the basin $B(K)$ of the statistical attractor $K$. Taking into account Equality (\ref{equationBsubEpsilonDeH}) which defines the set $B_{\epsilon}(K) \subset M$ for all $\epsilon >0$, we deduce that $B(K) = \bigcap _{\epsilon >0}B_{\epsilon}(K)$. From the statements (a) and (b) of Theorem \ref{theoremC}  and Definition \ref{definitionobservable}, we deduce that there exists $x \in B(K) \subset B_{\epsilon' }(K))$ and $\tilde \mu \in   {\mathcal L}^*(x) $ such that $|\int \psi \, d\mu - \int \psi \, d  \tilde \mu| < \epsilon$. Therefore, there exists a subsequence $\{\nu_{n_i}(x)\}_{i \geq 1}$ convergent to $\tilde \mu$ in the weak$^*$ topology of ${\mathcal M}$, where $\nu_n(x) := (1/n)\sum_{j=0}^{n-1} \delta_{f^j(x)}$. Thus:


 $$\int \psi \, d \tilde \mu = \lim _{i \rightarrow + \infty} \int \psi \, d \left (\frac{1}{n_i} \sum_{j= 0}^{n_i-1} \delta_{f^j(x)}   \right ) = \lim _{i \rightarrow + \infty} \frac{1}{n_i} \, \sum_{j= 0}^{n_i-1} \psi (f^j(x)) \geq $$


 $$ (1- \epsilon) \, \lim_{i \rightarrow + \infty} \frac{1}{n_j} \, \# \{0 \leq j \leq n_j -1: f^j(x) \in {\mathcal V}_{\epsilon'}    \}.$$


 Since $x \in B_{\epsilon'}(K)$, the limit    of the right term   in the above inequality, is equal to $1$ (recall Equality (\ref{equationBsubEpsilonDeH})).
 Therefore,
  $\int \psi \, d \tilde \mu  \geq  1 - \epsilon     $,
  and thus $\int \psi \, d   \mu    \geq 1- 2\epsilon.$
 Joining this latter result with Inequality (\ref{equation7}), we deduce that
  $ \mu(V_{\epsilon}) \geq  1-2 \epsilon  \ \forall \ \epsilon >0.  $
 Taking $\epsilon \rightarrow 0^+$   and taking into account that   the compact set $K $ is the decreasing intersection of the open sets $V_{\epsilon}$, we obtain:


 $$1 \geq \mu(K) = \lim _{\epsilon \rightarrow 0^+} \mu (V_{\epsilon})  = 1.$$


We have proved that $\mu(K) = 1$ for all $\mu \in {\mathcal O}_f(K)$, as wanted.

  Finally, it is left to prove that if $K$ is   minimal $\alpha$-obs with basin $B(K)$, then $K$ is the minimum compact set such that $\mu(K)= 1$ for all $\mu \in {\mathcal O}_f(K)$. Take a nonempty and compact set $K' \subset K$  such that $K  \setminus K' \neq \emptyset $. We shall prove that   $\mu(K') < 1$ for some $\mu \in {\mathcal O}_f(K)$.
  The minimality hypothesis on $K$ implies that  the set
 $B(K') $ (according to Definition \ref{definitionMilnorLike}), excludes a Lebesgue-positive set of points of $B(K)$. In other words,  $m(C) >0$, where $C:=  B (K) \setminus B (K') = \bigcup_{\epsilon >0} B(K) \setminus B_{\epsilon}(K') \subset B(K)$, with $B_{\epsilon}(K')$ satisfying Equality (\ref{equationBsubEpsilonDeH}).  Fix a point $x \in C$ and fix $\epsilon >0$ such that $x \not \in B_{\epsilon}(K')$. Choose a continuous real function $\psi \in C^0(M, [0,1])$ such that $\psi|_{K'} = 1$ and $\psi (y) = 0$ for all $y $ such that $\mbox{dist} (y, K') \geq \epsilon$. After Equality (\ref{equationBsubEpsilonDeH}), we obtain
 $ \ \  \lim \inf_{N \rightarrow + \infty} \omega_{\epsilon, N}(x, K', \epsilon) < 1  $   for all $x \in C$.
 In other words, there exists a sequence $n_i \rightarrow + \infty$ such that


 $$\lim _{i \rightarrow + \infty} \frac{1}{n_i} \# \{0 \leq j \leq n_i-1: \ \ \mbox{dist}(f^j(x), K') < \epsilon\} <1 \; .$$


  Therefore,


  $$\limsup_{i \rightarrow + \infty} \int \psi \, d \nu_{n_i}(x) := \limsup_{i \rightarrow + \infty} \int \psi \, d \left ( \frac{1}{n_i} \sum_{j=1}^{n_i} \delta_{f^j(x)}     \right)  
  =\limsup_{i \rightarrow + \infty} \frac{1}{n_i}\sum_{j=0}^{n_i -1} \psi  (f^j(x)) 
  $$


  $$ \leq \limsup _{i \rightarrow + \infty} \frac{1}{n_i} \# \{0 \leq j \leq n_i-1: \ \ \mbox{dist}(f^j(x), K') < \epsilon\}   \leq 1 - \epsilon. $$


  Taking  if necessary  a subsequence of $\{n_i\}_{i \geq 0}$ (which we still denote $\{n_i\}_{i \geq 0}$)  such that $\{\nu_{n_i}(x)\}_{i \geq 0}$ is convergent  in the weak$^*$ topology  to a probability measure $\mu  \in {\mathcal L}^*(x)$, we obtain:
   $\ \int \psi \, d \mu  = \lim _{i \rightarrow + \infty} \int \psi \, d \nu_{n_i}(x)  < 1\, . $
But, on the other hand, $\int \psi \, d \mu \geq \int_{K'} \psi \, d \mu = \mu(K')$. So $\mu (K') < 1$.

We have proved that  for all $x \in C$  there exists a measure
$\mu = \mu_x \in {\mathcal L} ^*(x)$ such that $\mu_x (K') < 1
$. Recall that $C \subset B(K)$ and  $m(C) >0$. After the
statement (a) of Theorem \ref{theoremC} (which we have already
proved), ${\mathcal L}^*(x) \subset {\mathcal O}_f(K)$ for
$m$-a.e. $x \in B(K)$. So, in particular, the above inclusion holds for
$m$-a.e. $x \in C$. We conclude that $\mu (K') < 1$ for
some $\mu \in {\mathcal O}_f(K)$  as wanted.
\hfill $\Box$

\vspace{.3cm}
Theorem \ref{corollaryMinimalAttractingB} states that, for any given forward invariant set $B$ with positive Lebesgue measure, there exists a   statistical attractor that is minimal attracting $B$. We will show  how this attractor can be constructed:

\begin{maintheorem} \label{CorollaryH}
Let $B \subset M$ be   a nonempty and forward  invariant set \em(i.e. $B \subset f^{-1}(B)$) \em
such that $m(B) >0$. Construct the set ${\mathcal O}_{f|B}$ of all the SRB-like measures of $f|_B$.
Then, the minimal compact set $K \subset M$ such that $\mu (K)= 1$
for all $\mu \in {\mathcal O}|_{f|B}$ is a  statistical attractor,
its   basin of attraction contains $B$ Lebesgue a.e., and $K$ is  minimal attracting $B$.
\end{maintheorem}

{\em Proof: }    After Theorem \ref{corollaryMinimalAttractingB} there
exists a statistical attractor $K'$ that is minimal attracting $B$, i.e. $B(K') \supset B \ m$-a.e.   It is enough to prove that $K' = K$.

Applying Lemma \ref{lemma3} we have that $B\subset B(K)$ Lebesgue a.e. Since $K'$  is   minimal  attracting $B$ (see Definition \ref{definitionMinimalAttractingB}), we deduce that $K' \subset K$.

Now, let us prove that $K \subset K'$.  Notice that the
set ${\mathcal O}_f(K')$ of probability measures satisfying Assertions (a) and (b) of
Theorem \ref{theoremC}, coincides with the set ${\mathcal
O}|_{f|B(K')}$ of Lemma \ref{lemma2}. After Assertion (c)
of Theorem \ref{theoremC} $\mu(K') = 1$ for all $\mu \in {\mathcal
O}_{f|B(K')}$. Since $B \subset B(K')$ then ${\mathcal O}_{f|B} \subset {\mathcal O}|_{f|B(K')}$. Therefore $\mu(K')= 1$ for all $\mu \in {\mathcal O}_{f|B}$. By hypothesis $K$ is the minimal compact subset of the space such that $\mu (K)= 1$ for all $\mu \in {\mathcal O}_{f|B}$.  We conclude that $K \subset K'$
as wanted.   \hfill $\Box$


\section{Proof of the Lebesgue-decomposition  of the space} \label{SectionDecomposition} \label{proofTheoremD}

In this section we    prove  Theorem \ref{TheoremD}.


{\em Proof: }

After Theorem \ref{theoremExistenceOfAttractors}, there exists an $\alpha$-minimal observable statistical attractor $K_1$. Then $m(B(K_1)) \geq \alpha$    Denote $\alpha_1 = \alpha$.


\noindent \textsc{1st. Step.}

\noindent Denote $r_1 = m(B(K_1)) \geq \alpha$. Either $r_1 = 1$ or   $1- \alpha <     r_1 < 1 $ or $\alpha \leq r_1 \leq 1 -\alpha $.

(1) In the first case,    Theorem \ref{TheoremD} becomes
trivially proved with $p= 1$.

(2) In the second case denote $\alpha_2 = 1 - m(B(K_1))$. Then
$0 < \alpha_2 < \alpha$. Consider the set $B := M \setminus
B(K_1)$. After Definition \ref{definitionMilnorLike} it is
standard to check that $f^{-1}(B(K_1)) = B(K_1)$. Therefore
$f^{-1}(B) = B$. After
Theorem \ref{corollaryMinimalAttractingB}    there exists a
  statistical attractor $K_2$ which is minimal attracting
$B$. In other words, $K_2$ is minimal such that $  B \subset B(K_2)   $ \ $m$-a.e. As $\alpha_2 =
m(B)$, the attractor $K_2$ is   $\alpha'_2$-obs minimal
for $f|B$.
Therefore, in this second case  Theorem \ref{TheoremD} is
proved with $p= 2$.

(3) In the third case, the set
$B := M \setminus B(K_1)$ has Lebesgue measure $m(B) = \alpha_2 \geq \alpha$.  After Theorem \ref{corollaryMinimalRestrictedM'} there exists a statistical attractor $K_2$ that is $\alpha$-obs. minimal for $f|B$. Now we go to the second step  by discussing again into three sub-cases:


\noindent \textsc{2nd. Step. }

Denote $r_2:= m(B(K_1) \cup B(K_2)) = m(B(K_1)) + m(B(K_2)\setminus B(K_1)) \geq 2 \alpha$.
 \begin{center}
 Either $r_2 = 1$, or $1 - \alpha < r_2 < 1 $ or $2 \alpha \leq r_2 \leq 1 - \alpha$.
  \end{center}

 (1') In the first case,  Theorem \ref{TheoremD} becomes trivially proved with $p= 2$.

 (2') In the second case,   Theorem \ref{TheoremD} becomes proved with $p= 3$, after the construction of a statistical attractor $K_3$ following the same arguments that were used in (2) to construct $K_2$.

 (3') In the third case, we can construct a minimal $\alpha$-obs statistical attractor $K_3$ for $f|_{  {\, M \setminus (B(K_1) \cup B(K_2))}}$, by  applying the same arguments that we used above in (3) to construct $K_2$. Now, we go to the following step, by discussing about the value of $r_3: = m (B(K_1) \cup B(K_2) \cup B(K_3))$.


 \noindent \textsc{Last  Step.}

 After   $p \geq 1 $ steps as above,   define the number


 $$r_p = m( \bigcup_{i=1}^p  B(K_i)) = \sum_{i=1}^p m(B(K_i)\setminus
 \cup_{j= 1}^{i-1}B(K_j)) \geq p \, \alpha.$$


 Since $r_p \leq 1$,
 the last step $p$  satisfies $p   \leq 1/\alpha$ and $1 - \alpha < r_p$.
 So, $p =  \mbox{Integer part}(1/\alpha)$. Therefore, in the last step $p$
 we always eventually drop in the cases   $ (1  )$ or   $(2)$.
 We conclude that   there exists a finite number  ($p$  or $p+1$)  of
 statistical attractors satisfying the statements (a), (b) and (c) of Theorem \ref{TheoremD}.
\hfill $\Box$


\section{Examples} \label{sectionAppendix}


 \begin{Example}   \textsc{(Hu-Young Diffeo)} \em
  \label{examplehuyoung}

Consider   the topologically transitive $C^{2}$ diffeomorphism $f$ studied in \cite{huyoung}: it acts in the 2-torus $\mathbf{T}^2$, and is obtained by an isotopy from a linear Anosov in such a way that the eigenvalues of $df$ at a fixed point $x_0$ are modified. Along the  contracting subspace  the eigenvalue is still   smaller than 1,  while in the eigendirection tangent to a topologically expansive (topologically unstable)
$C^1$ submanifold, the eigenvalue is weakened to become equal to 1. In \cite{huyoung}
it is proved, for such an $f$, that the sequence in Equality   (\ref{equation1})
of empirical probabilities  converges to $\delta_{x_0}$ in the space ${\mathcal M}$
of all the Borel probability measures (endowed with the weak$^*$-topology), for
Lebesgue a.e. $x \in \mathbf{T}^2$. In other words, $\delta_{x_0}$ is a physical measure,
the ergodic attractor is $K = \{x_0\}$ and its basin of attraction covers $\mathbf{T}^2$
  Lebesgue-a.e. Therefore, the frequency of visits to any arbitrarily
small neighbourhood of the fixed point is asymptotically equal to 1, for Lebesgue almost all the
initial states. The asymptotic frequency of visits to all the rest of the space is zero. In other words, $\{x_0\}$ is the unique $\alpha-$ observable minimal statistical attactor, for all $0 < \alpha \leq 1$.

Nevertheless, since $f$ is transitive, the unique (and thus minimal) Milnor's attractor
is the whole torus.
\end{Example}

 \begin{Example}  \textsc{(Bowen homeomorphism)}\label{exampleBowen} \em

 This example is attributed to Bowen   in \cite{takens} and \cite{godo},
 and  was also posed  in \cite{japon}.    Consider in a two dimensional
 manifold a non singular homeomorphism $f$ (namely $m(f^{-1}(B))= 0$
 if and only if $m(B)= 0$, where $m$ is the Lebesgue measure). Construct such an $f$ so that:

  (i) $f$ has three fixed points $x_1$, $x_2$ and $x_3$.

  (ii) When restricted to the union of  three small compact and pairwise disjoint  neighbourhoods $U_1$, $U_2$ and $U_3$ of $x_1$, $x_2$ and $x_3$ respectively, $f$ is a diffeo onto $f(U_1 \cup U_2 \cup U_3)$,  and   the fixed points $x_1$ and $x_2$ are hyperbolic saddles, while $x_3$ is a hyperbolic source.

    (iii)  $W^s_{1} \setminus \{x_1\} = W^u_2 \setminus \{x_2\}$,  $W^u_{1} \setminus \{x_1\} = W^s_2 \setminus \{x_2\}$.  We denote $W^{s,u}_{1,2}$ to   half-branches of the global one-dimensional stable and unstable manifolds of $x_{1,2}$ respectively. They are embedded topological arcs  of $C^1$ type in a neighbourhood of the saddles $x_{1,2}$. So $W^s_1 \cup W^s_2$ is a compact, simple and closed arc  which is the boundary of an open set $V$ homeomorphic to a 2-ball.

     (iv) The hyperbolic source $x_3$ is in $V$ and the orbits in $V\setminus\{x_3\}$ include $x_1$ and $x_2$ in their $\omega$-limit sets and have $\{x_3\}$ as   $\alpha$-limit set.

 Note that such a $C^0$ map $f$ can be constructed  for any previously specified values of the eigenvalues of $df$ at the two saddles $x_1$ and $x_2$, and after   an adequate choice of the values   $f(x)$ for $  x \in V \setminus (U_1 \cup U_2 \cup U_3)$.

 Let us consider the restricted dynamical system $f|_{\overline V}$. On the one hand and  from the topological viewpoint, all the orbits of $V\setminus \{x_3\}$ are attracted to (i.e. have $\omega$-limit set contained in)    the   boundary $\partial V$. This closed arc is the unique 1-obs. minimal Milnor's attractor of $f|_V$.   On the other hand, from the statistical viewpoint the behavior of the system is much more delicate (i.e. when looking the asymptotic behavior of the sequence of empirical probability measures defined in Equality (\ref{equation1})).  In fact,       necessarily one  and only one of the following properties (A), (B) or (C) holds, and  any of the three  is   realizable if the eigenvalues of $x_1$ and $x_2$ are adequately chosen and the $C^0$ map $f|_{\displaystyle{ V \setminus (U_1 \cup U_2)}}$ is  well constructed:


 \noindent (A) There exists a unique SRB  measure attracting
  $V \setminus \{x_3\}$ which is $\delta_{x_1}$ or ${\delta}_{x_2}$.
  In this case either $\{x_{1}\}$ or $\{x_2\}$ is an ergodic attractor,
  it is the unique statistical attractor  and the physical measure $\delta_{x_i}$
  is ergodic. We prove that this case is nonempty  (see the argument following the end of
   Example \ref{exampleS1} in this section).


 \noindent (B) There exists a unique SRB  measure $\mu$ attracting
 $V \setminus \{x_3\}$, which is $\mu = t \delta_{x_1} + (1- t) \delta_{x_2}$
 for some constant $0 < t < 1$. In this case (B), the   set $\{x_1, x_2\}$
    is an ergodic attractor for $f|_V$, it is the unique statistical attractor,
    and the physical measure $\mu$ is non ergodic. Moreover, for
    an adequate choice of the eigenvalues of $x_1$ and $x_2$
    one can obtain this property for any previously specified value of $t \in (0,1)$. The existence of examples in this case (B)
  is stated for instance in   Lemma Part (i) of page 457 in \cite{japon}.
   For the detailed   construction of an example  in this case,
   consider $\lambda= 1/\sigma$  in the Equalities  of Theorem 1
    of \cite{takens},  and  construct  $f$ such that it preserves
    area in both the disjoint compact neighbourhoods $U_1$ and $U_2$
    of the saddles, and is adequately $C^0$-chosen outside $U_1 \cup U_2$
    to have the two saddles in the omega-limit of all the orbits of $V \setminus \{x_3\}$.
    We note that, after a standard computation   that we sketch in the proof at the end
    of this section, one should construct $f$   contracting
         outside $U_1 \cup U_2$,
    so   the sequence (\ref{equation1}) is   convergent
    according to formulae of Theorem 1 of \cite{takens}
    (with the parameters $\lambda= 1/\sigma$ in those formulae).


 \noindent (C) There does not exist any physical measure, since for Lebesgue almost all
 the points   $x \in V   $,
  the limit set ${\mathcal L}^*(x)$ of the empirical distributions
  of Equality (\ref{equationL*x}) is a segment in the
  space $\mathcal{M}$ of probabilities. In other words,
  the sequence (\ref{equation1}) of empirical probabilities
  for $f|_{\overline V}$ does not converge for Lebesgue a.e. initial state.
  Thus, there does not exist  any ergodic attractor.
  The existence of $C^{2}$ examples in this case (C) is proved
   in \cite{takens} and \cite{godo} for which
   the set ${\mathcal O}_{f|\overline V}$ of SRB-like measures is
      a segment which is always properly contained in   $[\delta_{x_1}, \delta_{x_2}]
      \subset {\mathcal M}$.
    Nevertheless,  one can construct $f$  of $C^0$ class in $\overline V$
    such that the set of SRB-like measures for $f|_{\overline V}$
    is exactly the segment   $[\delta_{x_1}, \delta_{x_2}]$
    (see the remark at the end of this section).

  In this case (C),  there exist  uncountably many SRB-like measures
  for $f|_{\overline V}$ (after Theorem  1.7   of \cite{CE3}).
  All of them are   supported on $ \{x_1, x_2\}$, due to the Poincar\'{e}
  Recurrence Theorem.  After   Theorem \ref{CorollaryH} the set
  $\{x_1, x_2\}$ is a statistical attractor.    Besides,
  since the common minimal compact support of all the measures
  in ${\mathcal L}^*(x)$ is $ \{x_1, x_2\}$ for Lebesgue a.e.
  $x \in V$, this statistical attractor is the unique $\alpha$-obs minimal one,
  for all $0 <\alpha \leq 1$. In other words,  in this case (C)
  of example \ref{exampleBowen}, the unique $\alpha$-obs. minimal Milnor's
  attractor $\partial V$, and  the unique $\alpha$-obs. minimal statistical attractor,
  are different, while Pugh-Shub's ergodic attractors do not exist.

 \end{Example}

  Let us exhibit now an example that shows that if the purpose is to
  find the (always existing) statistical attractors of a $C^1$ map, even under
  the strong hypothesis of   uniform hyperbolicity, then
  the classic approach of searching  for the invariant probability measures
  that are absolutely continuous with respect to Lebesgue may become noneffective.
  Nevertheless, as a consequence of Theorem \ref{theoremC},
  there exists an optimal nonempty subset of   probability measures
  that describe   the statistics of Lebesgue almost all the orbits
  (see Definition \ref{definitionobservable}). In other words, for $C^1$ mappings that are not $C^1$ plus H\"{o}lder, the optimal probability measures are not
  necessarily   absolutely continuous with respect to Lebesgue.

  \begin{Example}   \textsc{(Campbell and Quas expanding maps)} \em

  \label{exampleS1}
  Let us consider a one-dimensional  $C^1$  map $f: S^1 \mapsto S^1$ on the circle $S^1$,
  which is expanding, i.e. $|f'(x)|>1$ for all $x \in S^1$.
  In Theorem 1 of \cite{campbellquas}, Campbell and Quas proved that $C^1$-generically
  there exists a unique physical measure $\mu$,
  that this measure $\mu$ is mutually singular with respect to Lebesgue,
  and that its basin of attraction covers Lebesgue almost all the points.
  This measure $\mu$  is supported on a compact subset $K \subset S^1$
  (non necessarily properly contained in $S^1$).
  So, this compact support $K$ is by definition an ergodic attractor.
  It is   the unique   statistical attractor and it is $\alpha$-obs minimal
  for all $0 <\alpha \leq 1$, since the basin of statistical attraction of $\mu$
  covers Lebesgue almost all the space.
  It is described by a single SRB-like measure which, in this case, is SRB.

\end{Example}

   \begin{Example} \textsc{(Quas  expanding maps)} \em
   \label{exampleQuas}

   In \cite{quas} Quas  gave a $C^1$-non generic example,  of an expanding
   map $f$ on the circle $S^1$ (which is $C^1$ but non $C^1$-plus-H\"{o}lder),
    exhibiting  a statistical behavior  that is rather opposite
   to that of the generic case of Campbell and Quas in Example \ref{exampleS1}.
   He constructed  such an
   $f$  preserving the Lebesgue measure $m$, but for which $m$ is non ergodic.
   So, after Birkhoff Theorem and after the Ergodic Decomposition Theorem,
   for $m$-almost every point $x \in S^1$  the set ${\mathcal L}^*(x)$
   (defined in Equality (\ref{equationL*x})) consists of  one ergodic
   component  of $m$. In this example, as in the general case, there exists
   a unique $1$-obs. minimal  statistical attractor (Theorem \ref{theoremExistenceOfAttractors}). But in this example,
   the set of all the SRB-like measures that describe completely the statistical behavior
   has more than one probability.  On purpose, all the SRB-like measures describing the statistics of the Ilyashenko's attractors of $C^1$ expanding maps on the circle, have also other good ergodic properties, from the viewpoint of the thermodynamic formalism. In fact, in \cite{CE4} it is proved that they all satisfy the Pesin's Entropy Formula.

  \end{Example}

\noindent {\bf  Proof of  the existence of Case (A) in Example \ref{exampleBowen}}

\em There exists an homeomorphism $f$ as in Example \em \ref{exampleBowen} \em such that


$${\mathcal{O}}|_{f|\overline V} = \{\delta_{x_2}\}.$$ \em


{\em Proof: }
  Choose   $f$ and the eigenvalues of $x_1$ and $x_2$ so that $f|_{U_1 \cup U_2}$ is $C^1$ and area conservative. Construct first   an area preserving map in a small neighbourhood of $\partial V$. Then   perturb $f$ near $\partial V$, in the $C^0$ topology, without changing $f|_{\displaystyle {\partial V \cup U_1 \cup U_2}}$,   to become hyper dissipative in a small neighbourhood    of a fundamental domain of $W^s_{x_2} \setminus (U_1 \cup U_2) $, and    not too much dissipative in  a small neighbourhood   of a fundamental domain of $W^s_{x_1}\setminus (U_1 \cup U_2) $. Precisely, construct this perturbation $f$ such that it satisfies the following property:

 At each return time  $n_i(x)$ to $U_2$ (of any orbit with initial state $x \in  V \setminus \{x_3\}$), and at each return time $n'_i(x)$ to $U_1$ such that $n_i(x) < n'_i(x) < n_{i+1}(x)$, consider the distances


\begin{equation}
 \label{equationsdi}
 d_i (x) := \mbox{dist}(f^{n_{i}}(x), W^s_{x_2}), \ \ \ d'_i(x):= \mbox{dist}(f^{n'_{i}}(x), W^s_{x_1})\end{equation}


\noindent Make $f$ to be $C^0$ in $V \setminus (U_1 \cup U_2)$ near  $\partial V$, so the above distances satisfy the  following inequalities:


 $$0 <d_{i+1}(x) < d'_i(x) ^{\displaystyle{-\log d'_i(x)}}, \ \ \ \ \frac{d'_i(x)}{3} \leq d'_{i+1}(x)  \leq \frac{d'_i(x)}{2} .$$


  At each visit $i$   to the set $U_2$, denote $N_ i(2)$ (depending on $x$) to the time length that the orbit of $x$ spends  inside $U_2$, and denote $N_i(1)$ to the time length that it spends inside $U_1$  after its $i$-th. visit to $U_1$. Up to a constant $k >0$, the number of iterates between the $i-$th. and the $(i+1)-$th. visit  to $U_2$ is $N_i(2) + N_i(1)$.  Besides, after a standard computation, we obtain $$N_i(2) \geq -c_2 \cdot \log {d_i} >   c_2 (-\log d'_i)^2, \ \ \ N_i(1) \leq - c_1 \cdot \log {d'_i},$$ for some positive constants $c_1$ and $c_2$. So, there exists $c >0$ such that

  $$ N_i(2) \geq c \, (N_i(1))^2 \ \ \forall \ \ i \geq 1.$$


   Consider the accumulated time average $\omega_n(U_1)$ inside $U_1$  of the finite piece of orbit from instant 0 up to instant $n \geq 1$ (namely, the relative frequency of staying in $U_1$).

    First, if $n$ is     exactly the end instant of the staying time  inside $U_1$ at the $m$-th. visit to $U_1$,  then $\omega_n(U_1) $ is computed as follows:


 $$\omega_n(U_1) = \frac{\sum_{i= 1}^{m } N_i(1)}{km + \sum_{i= 1}^{m} N_i(2) + \sum_{i= 1}^{m}N_i(1)} $$


 $$
 \frac{1}{\omega_n(U_1)} =  {\frac {km + \sum_{i= 1}^{m}  N_i(2) }{ \sum_{i=1}^m N_i(1)  }      + 1} \geq  {\frac {  \sum_{i= 1}^{m}  [N_i(1)]^2 }{ \sum_{i=1}^m N_i(1)  }     }.$$


 Since $N_i(1) \rightarrow + \infty$ when $i \rightarrow + \infty$, then  $1/\omega_n(U_1) \rightarrow + \infty$ when $m \rightarrow + \infty$ and so $\omega_n(U_1) \rightarrow 0$.

 Second, if $n$ is larger than the end instant $n'$ of the staying time inside $U_1$ at the $m$-th visit, but smaller than the next return time to $U_1$, then $\omega_n(U_1) = (n'/n) \, \omega_{n'}(U_1) \leq \omega_{n'}(U_1) \rightarrow 0$ when $m \rightarrow + \infty$.

 Third  and finally,   if $n$ is a  stoping time such that $f^n(x) \in U_1$ during the $m$-th. visit of the orbit to $N_1$, but $n   $ is smaller than the end instant $n'$ of the staying time $N_m$ inside $U_1$, then $0 < n'- n < N_m \leq c_1 (-\log d'_m) $. Since $d'_{i+1} \geq d'_i/3$ for all $i \geq 1$, we have $d'_m \geq (1/3^m) \, d'_1(x)$ for all $m \geq 1$. Thus, there exists a constant $K (x)>0$ such that $-\log d'_m \leq K (x) \cdot  m$ for all $m \geq 1$. This implies that $0 < n'-n < N_m \leq c'(x) \cdot m$  where $c'(x)= c_1 \, \cdot K (x)$. On the other hand $n \geq m$. Therefore


 $$\omega_n(U_1) = \frac{n'}{n} \  \omega_{n'}(U_1) = \omega_{n'}(U_1)\left (1 + \frac{n'-n}{n}\right ) \leq \omega_{n'}(U_1)\left (1 + c'(x)\right ) \rightarrow 0 $$


\noindent when $ m \rightarrow + \infty.$

  We have proved that   $\lim_{n \rightarrow + \infty} \omega_n(U_1) = 0$ for all $x \in V \setminus \{x_3\}$. Besides,     $\lim _n \omega_n(U_2) + \omega_n(U_1) = 1$. We deduce that $\lim_n\omega _n(U_2) = 1$. Since the argument above also holds (for the same $f$)  for any arbitrary   neighbourhood $U'_2$ of the saddle $x_2$, we obtain that the sequence (\ref{equation1}) converges to $\delta_2$, as wanted.  \hfill $\Box$


 \noindent {\bf Remark  about Case (C) of Example \ref{exampleBowen}}


\em There exists an homeomorphism   $f$ as in Example \em \ref{exampleBowen}, \em for which


$${\mathcal O}_{f|\overline V} =[\delta_{x_1}, \delta_{x_2}].$$
 \em


\noindent {\em Sketch of the proof. }
 Let us apply similar arguments to those of the
proof of case (A),     making $f$    hyper dissipative   near  $W^s(x_2) \setminus
(N_1 \cup N_2)$ but  also hyper dissipative near $W^s(x_1) \setminus (N_1 \cup N_2)$.
We deduce, adapting the computations in the proof of case (A), that the empirical
sequence (\ref{equation1}) will have at least two   convergent subsequences, one
converging to $\delta_{x_1}$ and the other to $\delta_{x_2}$. Fix a metric
$\mbox{dist}^*$ in the space ${\mathcal M}$ inducing the weak$^*$ topology.
After the convex-like property stated and proved in Theorem 2.1 of \cite{CE3},
for all $t \in [0,1]$ the limit set ${\mathcal L}^*(x)$ contains an invariant measure
$\mu_t(x)$ such that


\begin{equation}
\label{equalityultima}
\mbox{dist}^*(\mu_t(x), \delta_{x_1}) = t \mbox{dist}^*(\delta_{x_2}, \delta_{x_1}).
\end{equation}


From Poincar\'{e} Recurrence Theorem $\mu_t$ is supported on $\{x_1, x_2\}$,
so it is a convex combination of $\delta_{x_1}$ and $\delta_{x_2}$. But the unique such a
convex combination   satisfying Equality (\ref{equalityultima}), is
$\mu_t = t \delta_{x_1} + (1+t) \delta_{x_2}$, if the metric $\mbox{dist}^*$ is chosen
to depend linearly on $t$ for the measures in the segment $[\delta_{x_1}, \delta_{x_2}]$.
So ${\mathcal O}_{f|\overline V}= [\delta_{x_1}, \delta_{x_2}]$, as wanted.
\hfill $\Box$

\vspace{.2cm}

 \noindent {\bf Existence of  Case (B) of Example \ref{exampleBowen}}


\em For all $0 <t<1$ there exists an homeomorphism   $f$ as in Example \em \ref{exampleBowen}, \em for which


 $${\mathcal O}_{f|\overline V} =\{t\delta_{x_1} +
(1-t) \delta_{x_2}\}.$$
 \em


 {\em Proof: } 

  Applying similar arguments to those of the proof of case (A),   let us construct  a    weakly dissipative map $f$  near  $W^s(x_2) \setminus (N_1 \cup N_2)$, such that it is  also weakly dissipative near $W^s(x_1) \setminus (N_1 \cup N_2)$. Precisely, let us denote $d_i(x)$ and $d'_i(x)$ the distances defined in Equalities (\ref{equationsdi}) in the proof of case (A). We can   perturb a  map $f$ in the $C^0$ topology, in $V \setminus (U_1 \cap U_2)$, so that


   $$\frac{d'_i}{3} \leq d_{i +1} \leq \frac{d'_i}{2}, \ \ \ \ \frac{d _i}{3} \leq d'_{i } \leq \frac{d_i}{2} \ \ \ \forall \ i \geq 1.$$


  \noindent Recall that $f|_{U_1 \cup U_2}$ is area preserving. Adapting  standard computations obtained by applying Hartman-Grossman Theorem  inside the neighbourhoods $U_1$ and $U_2$ of the two saddles, we deduce that the staying times $N_i(1)$ and $N_i(2)$ (during the $i-th$ visit to $U_1$ and $U_2$ respectively) satisfy the following inequalities, for some positive constants $c $ and $ k'(x) $:


$$N_i(1) \leq c  \frac{\log d'_i} {\log \sigma_1} \leq k'(x) \frac{i}{\log \sigma_1} \leq N_i(1) +1   \ \ \forall \ i \geq 1  ,$$


$$N_i(2) \leq c  \frac{\log d _i} {\log \sigma_1} \leq k'(x) \frac{i}{\log \sigma_2} \leq N_i(2) +1   \ \ \forall \ i \geq 1  ,$$


\noindent where $\sigma_{1,2} >1$ are the expanding eigenvalues of the saddles $x_{1,2}$ respectively.

After similar computations to those in the proof of case (A), we deduce that the frequencies $\omega_n(U_1)$ and $\omega_n(U_2)$ of visits of the finite piece of orbit up to any stoping time $n \geq 1$, to the neighbourhoods $U_1$ and $U_2$ respectively, can be computed as follows:


$$\omega_n(U_{1,2}) \sim \frac{ \sum_{i=1}^m N_i(1,2)}{k m + \sum_{i= 1}^m N_i(2) + \sum_{i= 1}^m N_i(1)}   $$


\noindent where $k$ is a constant and $m$ is the number of visits to $U_2$ up to time $n$.
Thus,


$$\frac{1} {\omega_n(U_{1})} \sim 1 + \frac{km \log \sigma_{1}} { k'(x) \sum_{i=1}^m i} + \frac{\log \sigma_{1}}{\log \sigma_2} \rightarrow 1+ \frac{\log\sigma_1}{\log\sigma_2}$$


\noindent and analogously $\displaystyle{\frac{1}{\omega_n(U_2 )} \rightarrow 1+ \frac{\log\sigma_2}{\log\sigma_1}}$.

After checking that $1 = (1 + \log \sigma_1/ \log \sigma_2)^{-1} + (1 + \log \sigma_2/ \log \sigma_1)^{-1}$ we deduce that the empirical sequence (\ref{equation1}) will be   convergent to


$$t \delta_1 + (1-t) \delta_2 , \mbox{ where } t = \frac{1+ \log \sigma_2/ \log \sigma_1} {2 + \log \sigma_2/ \log \sigma_1 + \log \sigma_1/ \log \sigma_2}. $$


Since the eigenvalues $\sigma_{1,2} >1$ can be  arbitrarily chosen, the parameter $t$ can be equalled to any previously specified value in the open interval $(0,1)$.
 \hfill $\Box$

 \vspace{.5cm}

 \noindent{\bf  Acknowledgements:} We thank A. Gorodetsky and Y. Zhao for their valuable comments to an early version of this paper. We thank ANII and CSIC of the Universidad de la Rep\'{u}blica, Uruguay, for their partial financial support.



\label{lastpage}

\end{document}